\documentclass[10pt]{amsart}

\usepackage{hyperref}
\usepackage{amssymb, latexsym, amsmath, amsfonts}
\usepackage{amscd}
\usepackage{amsthm}

\providecommand{\cal}{\mathcal} 

\newcommand{\myamstitle}[3]{
  \title[#2]{#1}
  \author{Alexandru Scorpan}
  \address{Department of Mathematics, University of California, Berkeley\\ 970 Evans Hall, Berkeley, CA 94720}
  \email{scorpan@math.berkeley.edu}
  \urladdr{www.math.berkeley.edu/\textasciitilde scorpan}
  \date{#3}
  }

\theoremstyle{plain} 

\newtheorem{theorem}{Theorem}[section] 
\newtheorem{lemma}[theorem]{Lemma}
\newtheorem{proposition}[theorem]{Proposition}
\newtheorem{corollary}[theorem]{Corollary}

\theoremstyle{remark} 

\newtheorem{remark}[theorem]{Remark}

\newenvironment{Proof}[1]{\begin{proof}[Proof of #1]}{\end{proof}}

\newcommand{\Lemmaref}[1]{Lem\-ma~\ref{#1}}
\newcommand{\Propref}[1]{Pro\-po\-si\-tion~\ref{#1}} 
\newcommand{\Thmref}[1]{Theo\-rem~\ref{#1}}
\newcommand{\Rkref}[1]{Re\-mark~\ref{#1}}
\newcommand{\Corref}[1]{Co\-rol\-la\-ry~\ref{#1}}

\newcommand{\mraise}[2]{\raisebox{#1}{$#2$}}      

\newcommand{\sub}[1]{\raisebox{-2pt}{$\!_{#1}$} }  
\newcommand{\Sub}[1]{\raisebox{-2pt}{$\!_{\,#1}$} }

\newcommand{\adjustnabla}[1]{\sub{#1}} 

\newcommand{\goth}[1]{\mathfrak{#1}}              
\newcommand{\defemph}[1]{{\sffamily\slshape #1}}  

\renewcommand{\phi}{\varphi}

\newcommand{\ie}{{\it i.e.}~}

\newcommand{\st}{such that}
\newcommand{\Iff}{if and only if}

\newcommand{\LC}{Levi-Civit\`a}
\newcommand{\SW}{Sei\-berg--Witten}

\newcommand{\SD}{self-dual}
\newcommand{\ASD}{anti-self-dual}

\newcommand{\ac}{al\-most-com\-plex}
\newcommand{\riem}{Riemannian}
\newcommand{\herm}{Hermitian}
\newcommand{\kah}{K\"ahler}

\newcommand{\str}{struc\-ture}

\newcommand{\aR}{\mathbb{R}}
\newcommand{\Zi}{\mathbb{Z}}
\newcommand{\Ce}{\mathbb{C}}
\newcommand{\Ha}{\mathbb{H}} 

\newcommand{\Sph}[1]{{\mathbb S}^{#1}}  

\newcommand{\del}{\partial}           
\newcommand{\til}[1]{\widetilde{#1}}  
\renewcommand{\Bar}[1]{\overline{#1}} 
\newcommand{\rec}[1]{\tfrac{1}{#1}}   
\newcommand{\iso}{\approx}            
\newcommand{\tens}{\otimes}           

\newcommand{\maps}{\longmapsto}       
\newcommand{\longto}{\longrightarrow} 

\newcommand{\inner}[1]{\langle #1\rangle} 
\newcommand{\Inner}[1]{\bigl\langle #1\bigr\rangle} 
\newcommand{\Rinner}[1]{\inner{#1}\Sub{\aR}}
\newcommand{\RInner}[1]{\Inner{#1}\Sub{\aR}}

\newcommand{\norm}[1]{\left\|#1\right\|}  

\newcommand{\rest}[1]{|_{#1}}                          
\newcommand{\srest}[1]{\raisebox{-2pt}{$|_{#1}$}}      
\newcommand{\At}[1]{\raisebox{-2pt}{$\Big|_{#1}$}}     

\newcommand{\End}{\operatorname{End}} 
\newcommand{\Hom}{\operatorname{Hom}} 
\renewcommand{\Im}{\operatorname{Im}} 

\newcommand{\tr}{\operatorname{tr}}   

\newcommand{\T}[1]{T_{#1}}             

\newcommand{\TM}{\T{M}}

\newcommand{\uaR}{\underline{\aR}}     
\newcommand{\uCe}{\underline{\Ce}}     

\newcommand{\cont}{\,\lrcorner\,}      
\newcommand{\Nabla}[1]{\nabla\adjustnabla{#1}}

\newcommand{\lie}[1]{\mathfrak{#1}}    

\newcommand{\Ks}{K^{*}}                

\newcommand{\Ga}{{\mathcal G}}          
\newcommand{\Con}{{\mathcal C}onn}       

\newcommand{\Cl}{{\rm C}\ell}       
\newcommand{\Ccl}{\mathbb{C}\ell}   

\newcommand{\cli}{\,\mraise{3pt}{_{\bullet}}\, }     
\newcommand{\varcli}{\mbox{}_{^{\bullet}}} 

\newcommand{\Spinc}{\ensuremath{Spin^{\Ce}}}   

\newcommand{\pinc}{pin$^{\!{\bf C}}$} 
\newcommand{\spinc}{s\pinc}         

\newcommand{\W}{{\mathcal W}}           

\newcommand{\D}{{\mathcal D}}

\newcommand{\GW}{\Gamma(\W^{+})}    
\newcommand{\Wp}{\W^{+}}            
\newcommand{\Lp}{\Lambda^{+}}       

\newcommand{\hnabla}{\overrightarrow{\nabla}} 

\newcommand{\tnabla}{\til{\nabla}}  
\newcommand{\tNabla}[1]{\tnabla\adjustnabla{#1}}

\newcommand{\bnabla}{\smash{\Bar{\nabla}}}
\newcommand{\bNabla}[1]{\bnabla\adjustnabla{\!#1}}

\newcommand{\Tor}[1]{\operatorname{Tor}^{#1}} 

\newcommand{\tsum}{\operatorname{\mraise{1pt}{\sum}}} 

\newcommand{\ebasis}{\ensuremath{\{e_{1},e_{2},e_{3},e_{4}\}}} 

\newcommand{\abasis}{\ensuremath{\{a_{1},a_{2},a_{3},a_{4}\}}}

\newcommand{\PSO}{\ensuremath{PSO^{+}(\TM)}}

\newcommand{\B}{{\mathcal B}}
\newcommand{\Alt}{\operatorname{Alt}}
\newcommand{\Sym}{\operatorname{Sym}}
\newcommand{\Symo}{\operatorname{Sym}\sub{0}}

\begin{document}

\myamstitle{Spinors as automorphisms\\ of the tangent bundle}{Spinors as automorphisms}{April 26, 2002}

\keywords{spinor, four-manifold, \ac, symplectic, \kah}
\subjclass[2000]{Primary 53C27; Secondary 57N13, 32Q60, 53D05}

\begin{abstract}
We show that, 
on a $4$-manifold $M$ endowed with a \spinc-\str\ induced by 
an \ac\ \str, 
a \SD\ (= positive)
spi\-nor field $\phi\in\GW$ is the same as a 
bundle morphism $\phi:\TM\to\TM$ acting on the fiber by \SD\ 
conformal transformations, 
such that the Clifford multiplication 
is just the evaluation of $\phi$ on tangent vectors, 
and that the squaring map $\sigma:\Wp\to\Lp$ acts by 
pulling-back the fundamental form of the \ac\ \str.
We use this to detect \kah\ and symplectic structures.
\end{abstract}

\maketitle


\section{Introduction}

This paper is concerned with \spinc-\str s on $4$-manifolds, when the 
\spinc-\str\ is induced from an \ac\ \str.
The aim of the paper is two-fold. On the one hand, we will 
present a non-standard language for describing \SD\ (= positive) 
spinor fields
as automorphisms of the tangent bundle. 
On the other hand, using spinor fields to deform \herm\ \str s, 
we detect \kah\ \str s (\Thmref{thm-kahler}) 
and characterize symplectic \str s (\Corref{cor-more.symplectic}), 
expanding upon results from \cite{scorpan.harmonic}.

\bigskip

Let $M$ be an oriented $4$-manifold, endowed with a metric $g$. 
Using this metric, we can identify $g$-orthogonal \ac\ \str s with 
\SD\ $2$-forms of constant length $\sqrt{2}$. 
Thus, for example, for a non-zero \SD\ $2$-form $\alpha$, we can define its Chern class $c_1(\alpha)$ as the first Chern class of the associated \ac\ \str.

Choose an \ac\ \str\ $\omega$. It induces a standard \spinc-\str\ with 
spinor bundles denoted $\W^{\pm}$, with determinant bundle 
$\Ks=\det_{\Ce}\W^{\pm}=\det_{\Ce}(\TM,\omega)$, 
and with a Clifford multiplication 
denoted by $\TM\times\Wp\stackrel{\varcli}{\longto}\W^{-}$. It is a 
standard fact that $\W^{-}\iso(\TM,\omega)$ as complex bundles.
A choice of unitary connection $A$ on $\Ks$, together with the \LC\ connection 
$\nabla$ of $g$, induce unique unitary connections $\bnabla^{A}$ on 
$\Wp$ and $\tnabla^{A}$ on $\W^{-}$, related by 
\begin{equation}
 \label{eq-comp}
 \tnabla^{A}(v\cli\phi)=(\nabla v)\cli\phi+v\cli(\bnabla^{A}\phi)
\end{equation}
An important related object is 
the associated Dirac operator $\D^{A}:\GW\to\Gamma(\W^{-})$, defined by 
$\D^{A}\phi=\sum e_{k}\cli\bnabla^{A}_{e_{k}}\phi$
for any $g$-ortho\-nor\-mal frame $\{e_{k}\}$.
Another important object is the quadratic map $\sigma:\Wp\to\Lp$, 
which is famous for appearing in the \SW\ equations, but see also 
\cite[IV.10]{spingeom}.

A section of $\Wp$ will be called a \SD\ spinor field. 
(A more customary terminology would be ``positive spinor field''. We prefer 
to say ``\SD\ spinor field'', which is used in the classical 
paper \cite{AHS} and seems better suited to the peculiarities of 
dimension $4$ and to the phenomena described in this paper.)

\bigskip

A first aim of this paper is to
show that \SD\ spinor fields $\phi\in\GW$ can be identified 
with \defemph{\SD\ conformal transformations} 
\[ \phi:\TM\to\TM \]
\ie with bundle maps that act on each fiber 
by rotating a pair of orthogonal planes 
by a same angle, in directions compatible with the orientation of 
$M$, and then dilate/shrink the fiber by multiplying with a scalar. 
This bundle identification can be easily obtained from the 
starting steps of \cite{taubes.sw=gr, taubes.book}
(see \Rkref{rk-taubes}), 
but we strengthen it by noticing that the Clifford multiplication 
identifies with the evaluation, as 
\begin{equation}
 \label{eq-id}
 v\cli\phi\equiv\phi(v)
\end{equation}
and that the 
quadratic map $\sigma:\Wp\to\Lp$ can be described as giving the 
pull-back of the \ac\ form, as 
\[ \sigma(\phi)\equiv\rec{4}\,\phi^{*}\omega \] 
where $(\phi^{*}\omega)(v,w)=\omega(\phi v,\ \phi w)$.
This non-standard language is stated in \Thmref{thm-dictionary}.

Thus, one can use spinor fields to deform \ac\ \str s (with\-in a Chern 
class; see \ref{rk-lift}).

\bigskip

A second aim of this paper is to expand on the following result:
\begin{proposition}[\cite{scorpan.harmonic}]
\label{prop-old}
Consider a $4$-manifold $M$ endowed with a metric $g$ and with the 
\spinc-\str\ induced from an \ac\ \str\ $\omega$. 
Assume that $H^{2}(M;\Zi)$ has no $2$-torsion.
Then the equality $\alpha=\sigma(\phi)$ establishes
a bijection between:

{\bf A.}
The set of all \kah\ forms $\alpha$ with $c_{1}(\alpha)=c_{1}(\omega)$ 
and compatible with a metric \emph{scalar}-multiple of $g$; 
and the set of all gauge classes of pairs $(\phi,A)$ with 
$\phi$ nowhere-zero and $\bnabla^{A}\phi=0$.

{\bf B.}
The set of all symplectic forms $\alpha$ with 
$c_{1}(\alpha)=c_{1}(\omega)$ and compatible with a metric 
\emph{conformal} to $g$;
and the set of all gauge classes of pairs $(\phi,A)$ with 
$\phi$ nowhere-zero, $\D^{A}\phi=0$, 
and $\RInner{\bnabla^{A}\phi,\,i\phi}=0$.
\end{proposition}

(Here, ``gauge class'' means equivalent with respect to the 
action of the \defemph{gauge group} $\Ga=\{f:M\to\Sph{1} \}$ of $\Ks$.
It acts on $\Ks$ and $\W^{\pm}$ by scalar multiplication, 
and that induces an action on 
unitary connections on $\Ks$ and on sections of $\W^{\pm}$ (and thus on 
pairs $(\phi,A)$). See also \Rkref{rk-lift}.)

An immediate remark about \Propref{prop-old} is the lack of symmetry 
of (A) and (B): 
one statement deals with metrics scalar-multiple of $g$, 
the other with metrics conformal to $g$. 
Another is that the term 
$\Rinner{\bnabla^{A}\phi,\ i\phi}$, while formally clear, has a rather obscure intuitive meaning.
These remarks will be addressed as follows:

\bigskip

We will extend statement (A) from above 
to include all \kah\ forms compatible 
with metrics \emph{conformal} to $g$ (instead of merely 
scalar-multiple of $g$).
But the connections considered so far are not enough. 
We need a more general set of connections $\tnabla$ on $\W^{-}$
(called ``admissible connections'') 
that do not relate to any connections $A$ 
on $\Ks$, and do not correspond (via \eqref{eq-comp}) 
to connections $\bnabla$ on $\Wp$, 
but to connections $\bnabla$ on the larger bundle $\Hom(\TM,\TM)\supset\Wp$. 

Concretely, using the \herm\ isomorphism $\W^{-}\iso\TM$, 
a connection  $\tnabla$ on $\W^{-}$ is \defemph{admissible} 
if it is $\Ce$-linear for $\omega$ and $g$-metric. 
(The connections $\tnabla^{A}$ associated to connections $A$ on $\Ks$
need the extra condition 
$\tnabla^{A}\rest{\Lambda^{-}}=\nabla\rest{\Lambda^{-}}$; see 
\Lemmaref{lemma-spin.connection})
An admissible connection $\tnabla$ defines a connection $\bnabla$ on 
$\Hom(\TM,\TM)$ through 
$(\bnabla\phi)(v)=\tnabla(\phi v)-\phi(\nabla v)$.
The latter is simply a version of \eqref{eq-comp} read using \eqref{eq-id}, 
and is natural if we view $\phi:(\TM,\nabla)\to(\TM,\tnabla)$.

The extension of \ref{prop-old}.A is:%
\begin{theorem}
\label{thm-kahler}
Assume $H^{2}(M;\Zi)$ has no $2$-torsion.
The equality $\alpha=\sigma(\phi)$ establishes
a bijection between:
the set of all \kah\ forms $\alpha$ with $c_{1}(\alpha)=c_{1}(\omega)$ 
and compatible with a metric conformal to $g$; 
and the set of all gauge classes of pairs $(\phi,\tnabla)$ with 
$\tnabla$ admissible, $\phi$ nowhere-zero, and 
with 
\[ (\bNabla{X}\phi)(Y)=(\bNabla{Y}\phi)(X) \]
\end{theorem}
(Here again, the gauge group $\Ga$ acts by scalar multiplication on 
$\W^{-}=\TM$ and thus induces an action on connections $\tnabla$, and 
hence on pairs $(\phi,\tnabla)$. See also \Rkref{rk-lift}.)

\Thmref{thm-kahler} above could be read as a strong 
$4$-dimensional \spinc\ cousin of 
Proposition 9.10 from \cite[p.~340]{spingeom}. 
The latter states:
\emph{%
Let $M^{n}$ be endowed with a spin-\str, and let $\phi$ be a 
pure spinor. Then $\phi$ determines an integrable \ac\ \str\ \Iff\
$Y\cli\bNabla{X}\phi=X\cli\bNabla{Y}\phi$, for all $X, Y$ from the 
kernel of the map $\TM\tens\Ce\to\W$, $w\mapsto w\cli\phi$, 
and where $\bnabla$ is the unique spin-connection on the spinor bundle $\W$.}
Notice that this statement mentions only integrability, but not \kah.

\bigskip

This paper will also interpret the rather mysterious term 
$\Rinner{\bnabla^{A}\phi,\ i\phi}$ from \Propref{prop-old}.B.
It appeared there due to the formula
\begin{equation}
 \label{eq-D.formula}
 \norm{\phi}^{2} \D^{A}\phi = i\bigl( 2\,d^{*}\sigma(\phi) + 
  \Rinner{\bnabla^{A}\phi,\,i\phi} \bigr) \varcli\phi
\end{equation}
from \cite{scorpan.harmonic}.

We will show that the $1$-form $\Rinner{\bnabla^{A}\phi,\ i\phi}$
measures how close are the connections 
$\nabla$ and $\tnabla^{A}$ when compared through 
$\phi:(\TM,\nabla)\to(\TM,\tnabla^{A})$ (see \Lemmaref{lemma-minima}).
In standard terms, $\inner{\bnabla^{A}\phi,\,i\phi}$
measures how far is $\bnabla^{A}\phi$ from being minimal 
(when $A$ varies). 
(That is, $\inner{\bnabla^{A}\phi,\,i\phi}=0$ \Iff\ $\bnabla^{A}\phi$ has point-wise minimal length.)

\bigskip

Underlying the above discussion is the general comparison 
of $\nabla$ and $\tnabla^{A}$ via $\phi$.
It is governed by an analogue 
of a ``second fundamental form'':
\[ (\B^{A}\phi)\Sub{X}Y=\tnabla^{A}_{X}(\phi Y)-\phi(\Nabla{X}Y) \]
In standard terms, it is simply 
\[ (\B^{A}\phi)\Sub{X}Y=Y\cli\bnabla^{A}_{X}\phi \]
It naturally splits as:
\begin{equation}
 \label{eq-split}
 \B^{A}\phi
  =\Alt\B^{A}\phi + \Symo\B^{A}\phi + g \tens \rec{4}\D^{A}\phi
\end{equation}
where $\Alt\B$ is the alternating (skew-symmetric)
part, $\Symo\B$ is the trace\-less-symmetric part, 
and $g\tens\rec{4}\D^{A}\phi$ is the trace part of $\B$.

The skew-symmetric part $\Alt\B$
compares through $\phi$ the torsions of $\nabla$ and $\tnabla$, 
and underlies \Thmref{thm-kahler} above.
Namely, $\Alt\B\phi=0$ means $(\bNabla{X}\phi)(Y)=(\bNabla{Y}\phi)(X)$,
and in that case 
$\sigma(\phi)=\rec{4}\,\phi^{*}\omega$ is \kah\ for $\phi^{*}g$.
On the other hand, the symmetric part $\Sym\B$ can be understood as 
the linear extension of terms like $(\B^{A}\phi)\Sub{X}X$, 
which compare geodesics in the $X$-direction. 
The trace $\D^{A}\phi$ can thus be understood
as an average comparison of geodesics.
 
\bigskip

As mentioned, $\inner{\bnabla^{A}\phi,\,i\phi}$ vanishes exactly when 
$\B^{A}\phi$ is point-wise minimized.
Trying to minimize the other terms appearing in the splitting 
\eqref{eq-split} through variation of
$A$ among unitary connections on $\Ks$, we will show 
that the various minimizing connections must all sit on 
a \emph{single affine line} and are distanced at \emph{fixed ratios}
(see \Thmref{thm-minima.line}).
Therefore, if two such minimizing connections happen to coincide, 
then all of them must coincide.
 
The conditions from \ref{prop-old}.B can be read now as: both 
$\B^{A}\phi$ and $\D^{A}\phi$ are minimized by a same $A$. But then 
this $A$ must minimize all terms from \eqref{eq-split}. In particular 
it must minimize $\Alt\B^{A}\phi$, that is, minimize the torsion of 
the connection $\tnabla^{A}$ when viewed through $\phi$.

In conclusion, \Propref{prop-old}.B can be rephrased as:
\emph{The form $\sigma(\phi)=\rec{4}\,\phi^{*}\omega$ is symplectic \Iff\ there is a 
connection $A$ that simultaneously minimizes all components of 
$\B^{A}\phi=\bnabla^{A}\phi$.}

In particular, we have: 
\emph{The form $\phi^{*}\omega$ is symplectic 
\Iff\ there is a connection $\tnabla^{A}$ on $\TM$
which, viewed through $\phi:(\TM,\nabla)\to(\TM,\tnabla^{A})$, 
simultaneously is a closest match to $\nabla$ (\ie $\B^{A}\phi$ minimal) 
and has the torsion minimized (\ie $\Alt\B^{A}\phi$ minimal).} 
Or even:

\begin{corollary}
\label{cor-more.symplectic}
The manifold $(M,g,\omega)$ is almost-\kah\ \Iff\ there is a 
$\Ce$-linear $g$-metric connection $\tnabla$ on $\TM$ with 
$\tnabla\rest{\Lambda^{-}}=\nabla\rest{\Lambda^{-}}$, and that simultaneously is closest to $\nabla$ and has minimal torsion.
\end{corollary}

\noindent 
We hope this characterization might shed some new light on the nature of \ac\ \str s compatible with symplectic \str s.

\bigskip

In what follows, Section 2 will deal with proving the non-standard 
language for spinors, while Section 3 will detail its geometric 
ramifications.


\section{Dictionary}

Let $M$ be a closed oriented $4$-manifold, endowed with a fixed 
\riem\ metric $g$ and its \LC\ connection $\nabla$. 
The same notation ``$\nabla$'' 
will denote the connections induced by $\nabla$ on the tensor bundles of $M$.
Using the metric, we will systematically identify $\TM$ and $\TM^{*}$,
and their corresponding bundles of tensors, 
including $\Lambda(\TM)$ and $\Lambda(\TM^{*})$. 
Throughout the paper, a suddenly appearing ``$x$'' 
will simply mean a generic point of $M$.

An \defemph{\ac} \str\ is an automorphism $J:\TM\to\TM$ \st\ $J\circ 
J=-id$. All \ac\ \str s considered will be compatible with the chosen 
orientation of $M$, and will be $g$-orthogonal, \ie
$g(v,w)=g(Jv,Jw)$.
The metric allows us to identify all such \ac\ 
\str s $J$ with \SD\ $2$-forms $\omega\in\Gamma(\Lp)$ of constant length $\sqrt{2}$, 
\st\ $\omega(v,w)=g(Jv,w)$. 
(This is just a particular instance of the 
isomorphism $\Lambda^{2}(\aR^{4})\iso\lie{so}(4)$, identifying $2$-forms 
with skew-symmetric endomorphisms.)
The comp\-lex-line bundle 
$\Ks=\det_{\Ce}(\TM,\omega)$ is called the \defemph{anti-canonical\/} 
bundle of $(M,\omega)$. We denote by $c_{1}(\omega)$ the Chern class
$c_{1}(\Ks)=c_{1}(\TM,\omega)$.

A non-degenerate $2$-form $\omega$ with $d\omega=0$ is called 
\defemph{symplectic}. 
If a symplectic form is \SD\ and of constant length 
$\sqrt{2}$, then we call it a symplectic form \emph{compatible} 
with the metric $g$ 
(and in that case $(M,g,\omega)$ is an almost-\kah\ manifold).
A \SD\ $2$-form $\omega$ of length $\sqrt{2}$ \st\ $\nabla\omega=0$ will 
be called a \defemph{\kah\ form} compatible with $g$ 
(since $(M,g,\omega)$ is a \kah\ manifold).

\bigskip

The method we choose for proving \Thmref{thm-dictionary} below employs 
qua\-ternions, and is inspired by the exposition 
from \cite{akbulut}. 

Denote by $\Ha$ the division algebra of quaternions, and by $\Sph{3}$ 
its unit sphere. The choice of any isomorphism $\Ha\iso\aR^{4}$ 
that preserves orientation and inner-product allows us to identify 
$SO(4)=\Sph{3}\times\Sph{3}\,\big/\pm1$ acting on $\aR^{4}$ by 
\[ SO(4)\times\aR^{4}\longto\aR^{4} \quad:\quad
  [\xi_{+},\xi_{-}]\cdot v=\xi_{+}v\,\xi_{-}^{-1} \] 
If we further identify $\Ha$ with $\Ce^{2}$, 
through $z_{1}+z_{2}j\equiv(z_{1},z_{2})$, 
then we can identify $SU(2)=\Sph{3}$ acting by 
\[ SU(2)\times\Ce^{2}\longto\Ce^{2} \quad:\quad
  \xi\cdot v=v\xi^{-1} \]
The full unitary group can be identified as 
$U(2)=\Sph{1}\times\Sph{3}\,\big/\pm1$ acting on $\Ce^{2}$ by 
\[ U(2)\times\Ce^{2}\longto\Ce^{2} \quad:\quad
  [\lambda,\xi]\cdot v=\lambda v\xi^{-1} \]
The complex-spin group is 
$\Spinc(4)=\Sph{1}\times\Sph{3}\times\Sph{3}\,\big/\pm1$.

Since $M$ is oriented and endowed with a metric, its tangent bundle 
$\TM$ admits a defining cocycle with values in $SO(4)$. That
means that there is a covering of $M$ by open sets $\{U_{\gamma}\}$ 
and a collection of \emph{transition maps} 
$\{\tau_{\alpha\beta}:U_{\alpha}\cap U_{\beta}\to SO(4)\}$
\st\ the bundle $\TM\to M$ can be obtained by gluing trivial-bundle pieces 
$U_{\gamma}\times\aR^{4}\to U_{\gamma}$ through the identification of 
$(x,v_{\alpha})\in U_{\alpha}\times\aR^{4}$ with 
$(x,v_{\beta})\in U_{\beta}\times\aR^{4}$ when 
$v_{\alpha}=\tau_{\alpha\beta}(x) v_{\beta}$.

A \spinc-\str\ on $M$ is (the equivalence class of)
a lifting of the cocycle $\{\tau_{\alpha\beta}\}$ to a 
cocycle $\{\til{\tau}_{\alpha\beta}\}$ with values in $\Spinc(4)$, 
lifted via the natural map 
\[ \Spinc(4)\longto SO(4) \quad:\quad 
  [\lambda,\xi_{+},\xi_{-}]\maps[\xi_{+},\xi_{-}] \]
A choice of such a lifted cocycle $\{\til{\tau}_{\alpha\beta}\}$ 
induces, through the two maps 
\[ \Spinc(4)\longto U(2) \quad:\quad
  [\lambda,\xi_{+},\xi_{-}]\maps[\lambda,\xi_{\pm}] \]
defining $U(2)$-cocycles for two complex-plane bundles $\Wp$ and $\W^{-}$, 
called the bundles of \defemph{\SD} and \defemph{\ASD\ spinors}
(or ``positive'' and ``negative spinors'').
Through the map $\Spinc(4)\to\Sph{1}$, 
$[\lambda,\xi_{+},\xi_{-}] \maps \lambda^{2}$,
the cocycle $\{\til{\tau}_{\alpha\beta}\}$
also induces a defining cocycle for 
the complex-line bundle $L=\det_{\Ce}\W^{\pm}$,
which is called the \defemph{determinant} line-bundle of the \spinc-\str.
Every $4$-manifold admits at least one \spinc-\str.

The spinor bundles also come equipped with a Clifford 
multiplication $\TM\times\Wp\stackrel{\varcli}{\longto}\W^{-}$ (and 
its adjoint $\TM\times\W^{-}\to\Wp$), characterized by the property 
that 
\[ v\cli(v\cli\phi)=-\norm{v}^{2}\phi \]
This is in fact part of an action of 
the complexified Clifford algebra bundle $\Cl(\TM)\tens\uCe$ 
on $\Wp\oplus\W^{-}$ (see \cite[Ch.~II]{spingeom} for details). Since 
as vector bundles $\Cl(\TM)\iso\Lambda(\TM)$, there is an induced 
Clifford action of $\Lambda^{+}(\TM)$ on $\Wp$. 

It is known that, via Clifford multiplication, we have
$\End\sub{0}(\W^{+})\iso\Lambda^{+}\tens\Ce$ (where $\End\sub{0}$ 
denotes the set of \emph{traceless} $\Ce$-endomorphisms).
On the other hand, for every $\phi\in\Wp$, consider the endomorphism 
$\phi\tens\phi^{*}=\inner{\cdot,\phi}\phi$ of $\Wp$. Its traceless 
part is $\phi\tens\phi^{*}-\rec{2}\norm{\phi}^{2}id$. 
The latter corresponds to an element of $\Lambda^{+}\tens\Ce$, which 
turns out to be purely imaginary, \ie of the form $i\sigma(\phi)$ for 
some $\sigma(\phi)\in\Lambda^{+}$. This defines the \defemph{squaring map}
\[ \sigma:\Wp\to\Lp \]
Alternatively, $\sigma$ is uniquely characterized by its codomain and 
the property:
\[ \sigma(\phi)\cli\phi=-i\tfrac{\norm{\phi}^{2}}{2}\phi \]
This map is involved in the \SW\ equations, see for example
\cite{witten, donaldsonSW, morgan}.

\bigskip

If $(M,g)$ is endowed with a compatible \ac\ \str\ $\omega$, then the 
cocycle of $\TM$ can be reduced to a cocycle with 
values in $U(2)$. But there is a natural embedding 
$U(2)\subset\Spinc(4)$ given as 
\[ U(2)\longto\Spinc(4) \quad:\quad
  [\lambda,\xi]\maps[\lambda,\lambda,\xi] \]
Therefore the $U(2)$-cocycle 
of $\TM$ lifts to a canonical \spinc-\str\ 
$\{\til{\tau}_{\alpha\beta}\}$ associated with the \ac\ \str\ 
$\omega$.

Concretely, identify the model-fiber of $\TM$ with $\Ha$. 
If $\TM\rest{U}\iso U\times\Ha$ and $\TM\rest{U'}\iso U'\times\Ha$ 
are two bundle-charts,  
write $[\lambda,\xi]:U\cap U'\to\,\Sph{1}\times\Sph{3}\big/\pm1$ 
for the associated $U(2)$-transition map, which identifies 
$(x,v)\in U\times\Ha$ with $(x,v')\in U'\times\Ha$ when
$v'=\lambda(x)\,v\,\xi(x)^{-1}$.
The associated \spinc-\str\ will have corresponding $\Spinc(4)$-transition map 
$[\lambda,\lambda,\xi]:U\cap U'
  \to\,\Sph{1}\times\Sph{3}\times\Sph{3}\big/\pm1$.

The spinor bundle $\Wp$ will have model-fiber $\Ha$ and
corresponding $U(2)$-transi\-tion map
$[\lambda,\lambda]:U\cap U'\to\,\Sph{1}\times\Sph{3}\big/\pm1$, 
identifying $(x,w_{+})$ with $(x,w'_{+})$ when  
$w'_{+}=\lambda(x)\,w_{+}\lambda(x)^{-1}$.
The spinor bundle $\W^{-}$ will have model-fiber $\Ha$ and 
corresponding $U(2)$-transition map
$[\lambda,\xi]:U\cap U'\to\,\Sph{1}\times\Sph{3}\big/\pm1$
identifying $(x,w_{-})$ with $(x,w'_{-})$ when
$w'_{-}=\lambda(x)\,w_{-}\xi(x)^{-1}$.
The determinant line-bundle $L$ will have model-fiber $\Ce$ and 
corresponding $U(1)$-transition function
$\lambda^{2}:U\cap U'\to\Sph{1}$, 
identifying $(x,z)$ with $(x,z')$ when 
$z'=\lambda(x)^{2}z$.
(Note that, although the bundles $\W^{\pm}$ have model-fiber $\Ha$, 
after gluing them up with their respective cocycles, 
no global  quaternionic \str\ is preserved, 
only a complex \str.)

After inspecting the cocycles, immediate consequences are the well-known 
isomorphisms of \herm\ bundles $\W^{-}\iso\TM$ and $L\iso\Ks$.

\bigskip

Suppose now that $V$ is a $4$-dimensional vector space, 
endowed with an inner product and with an orientation.
Denote by $PSO(V)$ the group of orientation-preserving conformal 
transformations of $V$ (\ie real multiples of orthogonal 
transformations from $SO(V)$). 
Further, denote by $PSO^{+}(V)$ the group of 
\defemph{\SD\ conformal transformations}, that is 
those maps $V\to V$ which, with respect to \emph{orienting} orthonormal 
bases, are represented by matrices of the form:
\[\begin{bmatrix}
\phantom{-}r\cos\theta &  r\sin\theta & & \\
-r\sin\theta &  r\cos\theta & & \\
& & \phantom{-}r\cos\theta &  r\sin\theta\\
& & -r\sin\theta &  r\cos\theta
\end{bmatrix}\]
Note that the zero map $V\to 0$ is in $PSO^{+}(V)$, as is the 
identity $id:V\to V$.

Equivalently: Pick any vector-space isomorphism $V\iso\Ha$ preserving 
the inner product and the orientation. 
Then, a \SD\ conformal transformation is the same as the map $v\mapsto\xi v$
defined by multiplying on the left with a quaternion $\xi$.

A third description:
Any linear $\phi:V\to V$ induces a map 
$\phi:\Lambda^{2}(V)\to\Lambda^{2}(V)$ with $\phi(v\wedge 
w)=\phi(v)\wedge\phi(w)$. 
Since $V$ is $4$-dimensional, oriented, and has an inner product, 
we have a splitting $\Lambda^{2}(V)=\Lambda^{+}(V)\oplus\Lambda^{-}(V)$.
An automorphism $\phi:V\to V$ is a conformal transformation \Iff\ it 
preserves the splitting $\Lambda^{2}=\Lambda^{+}\oplus\Lambda^{-}$.
It is a \SD\ conformal transformation \Iff\ 
in addition it acts trivially on the 
\ASD\ part, that is, if $\phi\rest{\Lambda^{-}(V)}=id$.

Let $M$ be a closed oriented $4$-manifold, endowed with a fixed metric $g$. 
We can define the bundle $\PSO$ as the subbundle 
of $\Hom(\TM,\TM)$ containing the \SD\ conformal transformations of 
the fibers. 
It is a vector bundle of rank $4$.
As a subbundle of $\Hom(\TM,\TM)$, it comes equipped with a 
fiber-metric induced from $g$, as well as an
obvious evaluation map $\TM\times\PSO\stackrel{\rm ev}{\longto}\TM$, 
$(v,\phi)\maps\phi(v)$.

Suppose now that $M$ is endowed with some \ac\ \str\ $\omega$. 
Then $\TM$ becomes a complex bundle $(\TM,\omega)$.
A complex \str\ is induced on the bundle $\PSO$
simply by: $(i\cdot\phi)(v)=i\cdot\phi(v)$.
Denote the resulting \herm\ bundle by $(PSO^{+},\, \omega)$.

More, for any $\phi\in\PSO$, we can define the pull-back
$\phi^{*}\omega$ of the fundamental $2$-form $\omega\in\Lp(\TM^{*})$:
\[ (\phi^{*}\omega)(v,w)=\omega(\phi v,\phi w) \]
Since $\phi$ preserves $\Lp$, 
the pull-back $\phi^{*}\omega$ will be in $\Lp(\TM^{*})$
as well.

\begin{theorem}
\label{thm-dictionary}
Let $M$ be an oriented $4$-manifold endowed with a \riem\ metric and
a compatible \ac\ \str\ $\omega$, which induce a 
\spinc-\str\ with spinor bundles $\W^{+}$ and $\W^{-}$. 

{\bf A.}
We have the natural \herm\ bundle isomorphisms
\[ \W^{-}\iso(\TM,\,\omega) \qquad\qquad
  \Wp\iso(\PSO,\ \omega) \]
  
{\bf B.}
The Clifford multiplication $\TM\times\Wp\stackrel{\varcli}{\longto}\W^{-}$
identifies with the evaluation map 
$\TM\times\PSO\stackrel{\rm ev}{\longto}\TM$ as:
\[ v\cli\phi\equiv\phi(v) \]

{\bf C.}
The squaring map $\sigma:\Wp\to\Lp$ can be written:
\[ \sigma(\phi)=\rec{4}\,\phi^{*}\omega \]
\end{theorem}
\begin{remark}
Varying the \ac\ \str\ $\omega$ does not change the underlying real 
bundles of the spinor bundles, nor the Clifford multiplication map. 
It only changes the complex \str s that are laid on them.
\end{remark}
\begin{Proof}{\ref{thm-dictionary}.A}
We already proved $\W^{-}\iso\TM$. To show $\Wp\iso PSO^{+}$, we need only 
uncover the cocycle of the latter. 
Identify the model-fiber of $\TM$ with $\Ha$. Then 
a \SD\ conformal transformation of $\TM$ can be represented 
fiber-wise by left-multiplication with a quaternion. 
Pick a \SD\ conformal transformation $\phi:\TM\to\TM$, and
consider two bundle-charts $\TM\rest{U}\iso U\times\Ha$ and 
$\TM\rest{U'}\iso U'\times\Ha$ related by some $U(2)$-transition map 
$[\lambda,\xi]$. Let $x\in U\cap U'$.
If in the chart over $U$ the transformation $\phi$ is represented as
$v\mapsto hv$ for some $h:U\to\Ha$, 
while in the chart over $U'$ it is represented as 
$v'\mapsto h'v'$ for some $h:U'\to\Ha$, 
then, since the coordinate change is $v'=\lambda v \xi^{-1}$, 
we must have
$\lambda h v\xi^{-1}
  =h'\lambda v \xi^{-1}$ 
for all $v$. 
Therefore 
$h'=\lambda h \lambda^{-1}$.
In conclusion, $\PSO$ has the same transition functions as $\Wp$,
and hence these \herm\ complex bundles are isomorphic.
\end{Proof}

\begin{remark}
\label{rk-taubes}
The isomorphism $\W^{-}\iso\TM$ is well-known, usually written as 
$\W^{-}\iso\Lambda^{0,1}$. The isomorphism $\Wp\iso PSO^{+}$ follows 
immediately from the beginning of \cite{taubes.sw=gr, 
taubes.book}. For example: The Clifford action of $\omega$ splits 
$\Wp$ into $\pm1$-eigenbundles as $\Wp\iso\uCe\oplus\Ks$. But 
$\Lambda^{+}=\aR\omega\oplus\Ks$. Thus $\Wp\iso\uaR\oplus\Lp$.
Identifying $\Lambda^{+}$ with skew-symmetric endomorphisms and 
writing the trivial component as $\uaR=\aR\cdot id$ shows that 
$\Wp\iso PSO^{+}$ as real bundles. The complex structures follow as 
well.
\end{remark}
\begin{Proof}{\ref{thm-dictionary}.B}
The main task here is to concretely define the Clifford multiplication 
in such a manner that 
the identification $v\cli\phi=\phi(v)$ become obvious.

The action of the structure group $\Spinc(4)$ on $\Wp\oplus\W^{-}$ 
can be extended to an action of the full 
complexified Clifford algebra $\Cl(4)\tens\Ce$, 
in a way that respects the inclusion $\Spinc(4)\subset\Ccl(4)$. 
Globalizing we obtain the extension of
the action of the principal bundle $\Spinc(\TM)$ on $\Wp\oplus\W^{-}$ 
to an action of the algebra bundle $\Cl(\TM)\tens\uCe$.

But $\Cl(4)$ is isomorphic with the algebra $\Ha(2)$ 
of all $2\times 2$ quaternionic matrices, 
with $Spin(4)$ embedded as the group of all matrices
\[ \begin{bmatrix}\xi_{+}&\\&\xi_{-}\end{bmatrix} \]
with $\xi_{\pm}\in\Sph{3}$.
The inclusion of $\Spinc(4)=Spin(4)\times_{\Zi_{2}}\Sph{1}$ into 
$\Cl(4)\tens\Ce$ follows suit.

Identify the model-fiber of $\Wp\oplus\W^{-}$ with $\Ha\oplus\Ha$.
We define the Clifford multiplication in local quaternionic coordinates by
\[\begin{CD}
\Cl(\TM) &\ \tens\ & \uCe  \phantom{\At{}}
  &\qquad\times\qquad& \Wp\oplus\W^{-} 
  @>{\varcli}>> \Wp\oplus\W^{-}  \\ 
\begin{bmatrix}a&b\\c&d\end{bmatrix} &\tens& \lambda 
  && \begin{bmatrix}h^{+} & h^{-}\end{bmatrix}
  &
  & \lambda\cdot \begin{bmatrix}h^{+} & h^{-}\end{bmatrix} \cdot
     \begin{bmatrix} \Bar{a}&\Bar{c}\\ \Bar{b}&\Bar{d} \end{bmatrix}
\end{CD}\]
This local description  is compatible with the cocycles of $\TM$ and 
$\W^{\pm}$, and thus defines a global action of $\Cl(\TM)\tens\uCe$  
on $\Wp\oplus\W^{-}$.

The tangent bundle $\TM$ embeds into $\Cl(\TM)$ via a version of the 
(quaternionic) Pauli matrices.
Namely, for a local identification $\TM\rest{x}\iso\Ha$, we embed
$\TM\rest{x}$ into $\Cl(\TM)\rest{x}\iso\Ha(2)$ through:
\[ 1\mapsto\begin{bmatrix}0&-1\\ 1&\phantom{-}0\end{bmatrix},\  
   i\mapsto\begin{bmatrix}\phantom{-}0&-i\\ -i&\phantom{-}0\end{bmatrix},\ 
   j\mapsto\begin{bmatrix}\phantom{-}0&-j\\ -j&\phantom{-}0\end{bmatrix},\  
   k\mapsto\begin{bmatrix}\phantom{-}0&-k\\ -k&\phantom{-}0\end{bmatrix} \]
(These matrices generate $\Ha(2)$ and satisfy the Clifford relations 
$E_{k}\cdot E_{k}=-I\!d$ and $E_{j}\cdot E_{k}=-E_{k}\cdot E_{j}$. 
The above are essentially the negatives of the standard Pauli matrices.)
In short, the embedding is:
\[ v\maps\begin{bmatrix}&-v\\ \Bar{v}&\end{bmatrix} \]
Via this inclusion $\TM\subset\Cl(\TM)$, 
the tangent vectors will act on spinors as 
follows:
\[\begin{CD}
\TM &\ \times&\  \Wp @>{\varcli}>> \phantom{-}\ \W^{-}
  &\qquad\qquad&
  \TM &\ \times&\  \W^{-} @>{\varcli}>> \ \Wp\\
v && h^{+} && h^{+}v
  &&
  v && h^{-} && -h^{-}\,\Bar{v}
\end{CD}\]
In particular, $v\cli(v\cli\phi)=-\norm{v}^{2}\phi$, as needed.
If we read $\TM\times\Wp\longto\W^{-}$ through the isomorphisms
$\W^{-}\iso\TM$ and $\Wp\iso\PSO$, we obtain the map
$\TM\times\PSO\longto\TM$ given by $(v,h^{+})\maps h^{+}v$.
That is exactly the evaluation map read in coordinates.
\end{Proof}%

\begin{Proof}{\ref{thm-dictionary}.C}
We prove that $\sigma(\phi)=\rec{4}\,\phi^{*}\omega$.
As mentioned before, we use the metric to identify $2$-forms with 
skew-symmetric endomorphisms. Namely, the form 
$\gamma\in\Gamma(\Lambda^{2}(\TM^{*}))$ will correspond to the 
endomorphism of $\TM$ that satisfies
\[ \gamma(v,w)=g\bigl(\gamma(v),\,w\bigr) \]
(We will use the same letter for the $2$-form and for the morphism. 
The context or the specific number of arguments each takes should be enough to 
distinguish them.)
The essential ingredient of the proof is:
\begin{lemma}
\label{lemma-2form.action}
For any $\beta\in\Lambda^{+}(\TM)$ and $\phi\in\Wp$, we have
\[ \beta\cli\phi=-2\,\phi\circ\beta \]
\end{lemma}
Now, the squaring map $\sigma:\Wp\to\Lp(\TM)$
is characterized by:
\[ \sigma(\phi)\cli\phi=-i\tfrac{\norm{\phi}^{2}}{2}\,\phi \]
Using \ref{lemma-2form.action}, that translates to
$\phi\circ\sigma(\phi)=\tfrac{\norm{\phi}^{2}}{4}\,i\phi$, 
which can be written as:
\[ \phi\circ\sigma(\phi)=\tfrac{\norm{\phi}^{2}}{4}\,J\circ\phi \]
Therefore, thinking of $\phi$ as $\phi:\TM\to\TM$, 
and ignoring the factor $\norm{\phi}^{2}/4$ for a moment, 
we see that $\sigma(\phi)$ must determine the unique 
\ac\ \str\ on $\TM$ that will make 
$\phi:(\TM,\sigma(\phi))\to(\TM,J\/)$ be $\Ce$-linear.
When $\phi\srest{x}=0$, statement (C) is immediate, 
while in general we have:
\begin{align*}
\rec{4}\, \omega(\phi v, \phi w)
  &=\rec{4}\, g(J\phi v, \,\phi w)
   =\rec{\norm{\phi}^{2}}\,
      g\bigl( \tfrac{\norm{\phi}^{2}}{4}\, J\phi v, \,\phi w \bigr) \\
  &=\rec{\norm{\phi}^{2}}\,
     g\bigl( \phi \sigma(\phi) v, \,\phi w \bigr) 
   =g\bigl( \sigma(\phi)v, \,w \bigr) \\
  &=\sigma(\phi)(v,w)
\end{align*}
(where we used that $\phi$ is conformal, and hence 
$g(\phi x,\phi y)= \norm{\phi}^{2} g(x,y)$).
Thus $\rec{4}\,\phi^{*}\omega=\sigma(\phi)$.
\end{Proof}

\begin{remark}
\label{rk-lift}
Statement \ref{thm-dictionary}.C interprets
the square $\sigma(\phi)$ as a 
(weighted) pull-back of the \ac\ \str\ $\omega$ through 
$\phi:\TM\to(\TM,\omega)$.
On the one hand, 
it is worth quoting the following statement:

\emph{Let $\alpha$ be a \SD\ $2$-form, and assume that
$H^{2}(M\setminus\{\text{zeros of $\alpha$}\};\,\Zi)$ has no $2$-torsion.
There is a 
\SD\ spinor field $\phi\in\GW$ \st\ $\alpha=\phi^{*}\omega$
\Iff\ 
$c_{1}(\alpha\rest{\text{off zeros}})
 =c_{1}(\omega)\rest{\text{off zeros}}$.} 
\cite{scorpan.harmonic}

On the other hand, if we fix a suitable \SD\ $2$-form $\alpha$,
we can ask how unique is a spinor field $\phi$ \st\ $\alpha=\phi^{*}\omega$.
Such a $\phi:\TM\to\TM$ has the homothety ratio prescribed 
from $\norm{\phi}^{2}=2\sqrt{2}\norm{\phi^{*}\omega}$, and
must map the complex planes of $\alpha$ 
onto the complex planes of $\omega$.
Nonetheless, it has the freedom of rotating those planes.
Concretely, if $\phi^{*}\omega\rest{x}=\psi^{*}\omega\rest{x}$, 
then $\phi\rest{x}=e^{i\theta}\psi\rest{x}$ for some angle $\theta$.
This angular freedom can be factored out using the gauge group 
$\Ga=\{f:M\to\Sph{1} \}$.
\end{remark}%
\begin{Proof}{\Lemmaref{lemma-2form.action}}
We show that $\beta\cli\phi=2\,\phi\circ\beta$.
Through the Clifford action of $2$-forms on $\Wp\oplus\W^{-}$,
every $2$-form acts on a \SD\ spinor only through its \SD\ part
(the \ASD\ part acts trivially). 
For example, if 
$\phi\in\Wp\rest{x}$ and $\abasis$ is 
any orienting orthonormal basis in $\TM\rest{x}$, we have:
\[ (a_{1}\wedge a_{2})\cli\phi=(a_{3}\wedge a_{4})\cli\phi
  =\rec{2}(a_{1}\wedge a_{2}+a_{3}\wedge a_{4})\cli\phi \]

  Through the isomorphism $\Cl(\TM)\iso\Lambda(\TM)$, the algebra 
multiplication in $\Cl(\TM)$ can be expressed as 
$v\cdot\gamma = v\wedge\gamma - v\cont\gamma$,
for $v\in\TM$ and any $\gamma\in\Cl(\TM)$. 
Here $\cont$ is the interior product 
$v\cont\gamma=\gamma(v,\cdot,\ldots,\cdot)$
(identify $\Lambda(\TM)$ and $\Lambda(\TM^{*})$). 
If $v$ and $w$ are orthogonal, then 
$v\cont w=0$, and hence
$v\cdot w= v\wedge w$.
In particular,
$(a_{1}\wedge a_{2})\cli\phi=(a_{1}\cdot a_{2})\cli\phi$ and 
$(a_{3}\wedge a_{4})\cli\phi=(a_{3}\cdot a_{4})\cli\phi$.

Therefore we have 
$\psi=\rec{2}(a_{1}\wedge a_{2}+a_{3}\wedge a_{4})\cli\phi$
\Iff\ $\psi=(a_{1}\cdot a_{2})\cli\phi$ and 
$\psi=(a_{3}\cdot a_{4})\cli\phi$. 
Since $a_{k}\cdot a_{k}=-1$, that is the same as
$a_{1}\cli\psi=-a_{2}\cli\phi$ and 
$a_{3}\cli\psi=-a_{4}\cli\phi$. 
Using the identification $\Wp\iso\PSO$, that becomes:
\[ \psi=\rec{2}(a_{1}\wedge a_{2}+a_{3}\wedge a_{4})\cli\phi
  \quad\iff\quad
  \psi(a_{1})=-\phi(a_{2}) 
  \ \ \&\ \  
  \psi(a_{3})=-\phi(a_{4}) \]

The element
$\alpha=a_{1}\wedge a_{2}+a_{3}\wedge a_{4}$ of $\Lp\rest{x}$
corresponds to the rotation
$\alpha:\TM\rest{x}\to\TM\rest{x}$ acting by $a_{1}\mapsto a_{2}$, 
$a_{2}\mapsto-a_{1}$, and $a_{3}\mapsto a_{4}$, $a_{4}\mapsto-a_{3}$.
The composition $\phi\circ\alpha$ will still be \SD, and will act  by: 
$(\phi\circ\alpha)(a_{1})=\phi(a_{2})$ and 
$(\phi\circ\alpha)(a_{3})=\phi(a_{4})$. 
Comparing with $\psi$, we conclude: 
\[ \alpha\cli\phi=-2\,\phi\circ\alpha \]

To get \Lemmaref{lemma-2form.action}, we need only remark that
for any $\beta\in\Lambda^{+}\rest{x}$, 
there is always a suitable basis 
$\{b_{1}, b_{2}, b_{3}, b_{4}\}$ in $\TM\rest{x}$ so 
that $\beta=r(b_{1}\wedge b_{2}+b_{3}\wedge b_{4})$, with 
$r=\rec{\sqrt{2}}\norm{\beta}$, and thus the above applies.
\end{Proof}


\section{Geometry}

Using the language of \Thmref{thm-dictionary}, we view a \SD\ spinor 
field as a bundle morphism $\phi:(\TM,g,\nabla)\to(\TM,g,\omega)$.
The connections $\tnabla$ on $\W^{-}$ are seen as connections on the target of 
such $\phi$.

Choose any unitary connection $A$ on $\Ks$. 
This $A$ can be combined with the \LC\ connection $\nabla$ to induce 
unitary connections $\bnabla^{A}$ on $\Wp$ and $\tnabla^{A}$ on $\W^{-}$.
That is done by lifting $A$ through the map $\Spinc(4)\to U(1)$, 
combining with the lift of $\nabla$ through $\Spinc(4)\to SO(4)$, and 
then projecting the combination through the two maps 
$\Spinc(4)\to U(2)$.
(For details, see \cite{spingeom, morgan}, or the proof of 
\Lemmaref{lemma-spin.connection}, at the end of this paper.)
We call such connections on $\W^{\pm}$ \defemph{spinorial} connections.

Viewing a spinorial connection $\tnabla$ on $\W^{-}$ as a connection on $\TM$ 
(via $\W^{-}\iso\TM$), it is obvious that $\tnabla$ is 
$\Ce$-linear for $\omega$ and $g$-metric. In fact:
\begin{lemma}%
\label{lemma-spin.connection}
Let $D$ be any connection on $\TM$. 
There is a connection $A$ on $\Ks$ \st\ $D=\tnabla^{A}$ 
(via the identification $\TM\iso\W^{-}$) \Iff\ $D$ is $g$-metric, $\Ce$-linear for 
$\omega$, and the induced connection $D\rest{\Lambda^{-}}$ 
on $\Lambda^{-}(\TM^{*})$ coincides with the one induced by the \LC\ 
connection $\nabla$.
\end{lemma}
\begin{remark}
\label{rk-find.A}
Given $D$ as above, one can find the suitable $A$ as 
follows: The fact that $D$ is $\Ce$-linear for 
$\omega$ can also be written as $D\rest{\Lambda^{+}}\,\omega=0$. 
Since $\Lambda^{+}=\aR\omega\oplus\Ks$, that implies that 
$D\rest{\Lambda^{+}}=\del\oplus A$ for some unitary connection $A$ 
on $\Ks$ (where $\del$ denotes the trivial connection on the 
trivialized bundle $\aR\omega$). 
If the conditions from \ref{lemma-spin.connection} are met, 
then $D=\tnabla^{A}$.
\end{remark}
\begin{remark}
The fact that connections $D$ with 
$D\rest{\Lambda^{-}}=\nabla\rest{\Lambda^{-}}$ play such an important 
r\^ole in spin-geometry could be justified in the language of 
\ref{thm-dictionary} as follows: A bundle morphism $\phi:\TM\to\TM$ 
is a \SD\ spinor field \Iff\ the induced map on $\Lambda^{2}$ 
preserves the splitting $\Lambda^{2}=\Lambda^{+}\oplus\Lambda^{-}$ and 
acts on $\Lambda^{-}$ as the identity. It is thus natural that the most
natural connections on the target of $\phi:(\TM,\nabla)\to\TM$
coincide with $\nabla$ on $\Lambda^{-}$. 
\end{remark}

\Lemmaref{lemma-spin.connection} suggests a non-standard approach: 
Instead of choosing $A$ and building $\tnabla^{A}$ and $\bnabla^{A}$, 
one could start with a connection $\tnabla$ on $\TM$ which is 
$g$-metric, $\Ce$-linear, and has 
$\tnabla\rest{\Lambda^{-}}=\nabla\rest{\Lambda^{-}}$. 
We know that $\tnabla=\tnabla^{A}$ for some $A$, but we do not 
determine $A$. 
Instead, seeing spinor fields as $\phi:(\TM,\nabla)\to(\TM,\tnabla)$ 
and using \eqref{eq-comp} in the form 
\begin{equation}
 \label{eq-my.comp}
 (\bnabla\phi)v=\tnabla(\phi v)-\phi(\nabla v) 
\end{equation}
we can define a connection $\bnabla$ on the whole $\Hom(\TM,\TM)$. 
Because of \ref{lemma-spin.connection}, we know that this $\bnabla$ 
must preserve the subbundle $\Wp\subset\Hom(\TM,\TM)$. 
The restriction of $\bnabla$ to $\Wp$ is then exactly $\bnabla^{A}$.

The \herm\ identification $\W^{-}\iso(\TM,g,\omega)$ 
further tempts one to consider 
connections $\tnabla$ on $\W^{-}$ that are merely $g$-metric and 
$\Ce$-linear for $\omega$. We call them \defemph{admissible} 
connections on $\W^{-}$.

Admissible connections no longer correspond to a connection $\bnabla$ on 
$\Wp$. None\-theless, using \eqref{eq-my.comp}, we can still define a 
connection $\bnabla$ on $\Hom(\TM,\TM)$, but it will no 
longer preserve the subbundle $\Wp$. 
The meaning of $\bnabla$ (be it spinorial or not), as read from 
\eqref{eq-my.comp}, is to compare the connections $\nabla$ and 
$\tnabla$ through $\phi:(\TM,g,\nabla)\to(\TM,g,\omega,\tnabla)$.
For easier manipulation, we introduce the notation 
\[ \bigl(\B^{\tnabla}\phi\bigr)\Sub{X}Y
  =(\bNabla{X}\phi)(Y)=\tNabla{X}(\phi Y)-\phi(\Nabla{X}Y) \]
which can be though of as the analogue of a ``second fundamental form'' 
for $\phi:\TM\to\TM$.

We already encountered the pull-back $\phi^{*}\omega=4\,\sigma(\phi)$ 
of the fundamental form. 
We can also pull-back the metric $g$ to $\phi^{*}g$, given by  
$(\phi^{*}g)(v,w)=g(\phi v,\phi w)$. Since $\phi$ is conformal, it 
is simply $\phi^{*}g=\norm{\phi}^{2}g$. We can as well pull-back the 
unitary connection $\tnabla$ to $\phi^{*}\tnabla$, defined by 
$\bigl(\phi^{*}\tnabla\bigr)\Sub{X}Y=\phi^{-1}\tNabla{X}(\phi  Y)$, 
where $\phi^{-1}$ is the inverse of $\phi$ (defined only off the 
zeros of $\phi$). Off the zeros of $\phi$, the comparison form 
$\B^{\tnabla}\phi$ can then be written 
$(\B^{\tnabla}\phi )\Sub{X}Y
 =\phi\bigl( (\phi^{*}\tnabla)\Sub{X}Y-\Nabla{X}Y \bigl)$.  

The zeros of $\phi$ create singularities, but, nonetheless, 
since $g$, $\omega$, and $\tnabla$ were compatible, 
so will their pull-backs. 
Namely: the form $\phi^{*}\omega$ is \SD\ 
and has constant length $\sqrt{2}$ for the (singular) metric $\phi^{*}g$, and 
thus corresponds to a (singular) \ac\ \str; 
the (singular) connection $\phi^{*}\tnabla$ is 
\herm, \ie\ it is $\phi^{*}g$-metric and $\Ce$-linear for 
$\phi^{*}\omega$. 

\begin{remark}
The zeros of spinors are hard to control. Even in the case of harmonic 
spinors (spinors fields $\phi$ for which there is an $A$ \st\ 
$\D^{A}\phi=0$), the zero-set is a countable $2$-rectifiable set, and 
thus has Hausdorff dimension as high as $2$ (see \cite{baer}).
\end{remark}

In general, the connection $\phi^{*}\tnabla$ will have 
torsion. But if, for example, $\B^{\tnabla}\phi= 0$, then $\bnabla\phi= 0$, 
and so $\phi$ has constant length. 
If further $\phi$ is non-trivial, then 
$\phi^{*}\tnabla$ is well-defined on all $M$,
and we have $\phi^{*}\tnabla=\nabla$.
In this case $(M,\,\phi^{*}g,\,\phi^{*}\omega)$ is a \kah\ manifold.
This suggests that the comparison form $\B\phi$ has control over
the geometry of the deformed structure $\phi^{*}\omega$.

\bigskip

If $V$ is a vector-space endowed with an inner product $q$, we have
the standard splitting 
$V\tens V=\Lambda^{2}(V)\oplus {\rm S}^{2}_{0}(V)\oplus \aR q$.
The comparison form $\B^{\tnabla}\phi$ is tensorial 
--- a section in $\TM^{*}\tens \TM^{*}\tens \TM$.
Applying the above splitting on the $\TM^{*}\tens\TM^{*}$-factor, we 
get the splitting
\[ \B^{\tnabla}\phi
  =\Alt\B^{\tnabla}\phi + \Symo\B^{\tnabla}\phi 
  + g \tens \rec{4}\tr(\B^{\tnabla}\phi) \]
which breaks $\B^{\tnabla}\phi$ 
into its skew-symmetric, traceless-symmetric, and trace parts.
Concretely, these are defined as:
\begin{align*}
& (\Alt\B^{\tnabla}\phi)\Sub{X}Y
  =\rec{2}\bigl( (\B^{\tnabla}\phi)\Sub{X}Y
    - (\B^{\tnabla}\phi)\Sub{Y}X \bigr) \\
& (\Sym\B^{\tnabla}\phi)\Sub{X}Y
  =\rec{2}\bigl( (\B^{\tnabla}\phi)\Sub{X}Y 
    + (\B^{\tnabla}\phi)\Sub{Y}X \bigr) \\
& \Symo\B^{\tnabla}\phi
  =\Sym\B^{\tnabla}\phi - g\tens \rec{4}\tr(\B^{\tnabla}\phi)\\
& \tr\B^{\tnabla}\phi
  =\Inner{\B^{\tnabla}\phi,\,g}=\tsum \,(\B^{\tnabla}\phi)\Sub{e_{k}}e_{k} 
\end{align*}
where $\ebasis$ is any $g$-orthonormal basis in 
$\TM\rest{x}$.
When $\tnabla=\tnabla^{A}$ is spinorial,
we denote $\B^{\tnabla}$ by $\B^{A}$, as expected. 
Notice that in that case
$\tr(\B^{A}\phi)$ coincides with $\D^{A}\phi$, the Dirac operator.

\bigskip

The skew-symmetric component $\Alt\B\phi$ is easy 
to interpret. Indeed, off the zeros of $\phi$, the following 
calculation holds:
\begin{align*}
( &\Alt\B^{\tnabla})\Sub{X}Y 
 =\rec{2}\bigl( \tNabla{X}(\phi Y)-\tNabla{Y}(\phi X)
  -\phi(\Nabla{X}Y)+\phi(\Nabla{Y}X) \bigr)\\
&=\rec{2}\bigl( \phi\phi^{-1}\tNabla{X}\phi Y 
  -\phi\phi^{-1}\tNabla{Y}\phi X
  -\phi[X,Y] -\phi(\Nabla{X}Y-\Nabla{Y}X-[X,Y]) \bigr)\\ 
&=\rec{2}\,
  \phi \bigl( \Tor{\,\phi^{*}\tnabla}(X,Y)-\Tor{\nabla}(X,Y) \bigl)
\end{align*}
where $\Tor{D}(X,Y)=D_{X}Y-D_{Y}X-[X,Y]$ denotes the torsion of a 
connection $D$.
Thus $\Alt\B^{\tnabla}\phi$ is the torsion-comparing component of 
$\B^{\tnabla}\phi$. 
Since $\nabla$ has no torsion, we simply have:
\[ (\Alt\B^{\tnabla}\phi)\Sub{X}Y
  =\rec{2}\,\phi\bigl( \Tor{\,\phi^{*}\tnabla}(X,Y) \bigr) \]
We are now ready to prove \Thmref{thm-kahler}, stated as:

\emph{%
The equality $\alpha=\phi^{*}\omega$ establishes
a bijection between:
the set of all \kah\ forms $\alpha$ with $c_{1}(\alpha)=c_{1}(\omega)$ 
and compatible with a metric conformal to $g$; 
and the set of all gauge classes of pairs $(\phi,\tnabla)$ with 
$\tnabla$ admissible, $\phi$ nowhere-zero, and 
with $\Alt\B^{\tnabla}\phi=0$.}

\begin{Proof}{\Thmref{thm-kahler}}
We have $\Alt\B^{\tnabla}\phi=0$ \Iff\ 
$\phi^{*}\tnabla$ is torsion-free.
But in that case $\phi^{*}\tnabla$ is the \LC\ connection of 
$\phi^{*}g$, and since it is also $\Ce$-linear 
for $\phi^{*}\omega$, we conclude that
the form $\phi^{*}\omega$ is \kah\ for the 
metric $\phi^{*}g$.

Conversely, assume $M$ admits a \kah\ form $\alpha$ 
for a metric $g'$ conformal to $g$, 
and has the same Chern class as $\omega$. 
The latter implies (by \Rkref{rk-lift}) that there is 
a nowhere-zero spinor field $\phi$ \st\ $\phi^{*}\omega=\alpha$. 
Since $\alpha$ has length $\sqrt{2}$ for $g'$, we deduce that we must 
also have $\phi^{*}g=g'$. 
If $\nabla'$ denotes the \LC\ connection of 
$g'$, then $\nabla'\alpha=0$.
Define $\tnabla=\phi_{*}\nabla'$, where 
$(\phi_{*}\nabla')\sub{X}Y=\phi\nabla'\sub{\!\!X}(\phi^{-1}Y)$.
Then $\tnabla$ is an admissible connection on $\TM$ 
(it is $g$-metric and $\Ce$-linear for $\omega$), 
and obviously has $\phi^{*}\tnabla=\nabla'$.
The latter being torsion-free, we must have $\Alt\B^{\tnabla}\phi=0$.
The gauge-invariance part follows easily.
\end{Proof}%

\begin{remark}
Unless $\phi$ has constant-length, or, equivalently, 
unless the metric $g'=\phi^{*}g$ is a scalar-multiple of $g$
(as in \ref{prop-old}), 
the connection $\tnabla$ is \emph{not} spinorial.
\end{remark}

Combining the bijections from 
\Thmref{thm-kahler} and \Propref{prop-old}
yields a coarse constraint that $\Alt\B$ imposes on the whole $\B$:

\begin{lemma}%
Assume that, for a spinor field $\phi$, we have $\Alt\B^{A}\phi= 0$.
If $\phi$ is non-trivial and has constant length, then $\B^{A}\phi= 0$.
\end{lemma}%

\bigskip

We now restrict our attention to the standard spinorial connections 
$\bnabla^{A}$ and $\tnabla^{A}$. Varying these means varying $A$. 
Denote by $\Con(\Ks)$ the set of all unitary connections on $\Ks$. 
Any two connections from $\Con(\Ks)$ differ by a 
global imaginary $1$-form.
When $A$ varies in $\Con(\Ks)$, that is 
when $A$ changes to $A+2i\theta$ for some $\theta\in\Gamma(\TM^{*})$, 
then $\tnabla^{A}$ changes to 
$\tnabla^{A+2i\theta}=\tnabla^{A}+i\theta$. 
The associated tensors change as:
\begin{equation}
\label{eq-change}
\begin{split}
 & (\B^{A+2i\theta}\phi)\Sub{X}Y
  =(\B^{A}\phi)\Sub{X}Y
  +i\theta(X)\phi(Y) \\
 & (\Alt\B^{A+2i\theta}\phi)\Sub{X}Y
  =(\Alt\B^{A}\phi)\Sub{X}Y
  +\tfrac{i}{2}\theta(X)\phi(Y)-\tfrac{i}{2}\theta(Y)\phi(X) \\
 & (\Sym\B^{A+2i\theta}\phi)\Sub{X}Y
  =(\Sym\B^{A}\phi)\Sub{X}Y
  +\tfrac{i}{2}\theta(X)\phi(Y)+\tfrac{i}{2}\theta(Y)\phi(X)\\
 & \D^{A+2i\theta}\phi
  =\D^{A}\phi+\tsum\, i\theta(e_{k})\phi(e_{k})
   =\D^{A}\phi +i\theta\cli\phi
 \end{split}
\end{equation}
A simple consequence of these formulae is:

\emph{%
If, for two connections $A$ and $A'$, we have 
$\B^{A}\phi=\B^{A'}\phi$, or 
$\Alt\B^{A}\phi=\Alt\B^{A'}\phi$, or  
$\Sym\B^{A}\phi=\Sym\B^{A'}\phi$, or 
$\D^{A}\phi=\D^{A'}\phi$, 
then we must have $A=A'$ on the support of $\phi$.
}

The same formulae also yield:

\begin{theorem}%
\label{thm-minima.line}
Let $(\B^{A}\phi)\Sub{X}Y=Y\cli\bnabla^{A}_{X}\phi$ and vary $A$.
For every nowhere-zero spinor field  $\phi$, there is a unique connection $A_{0}$ 
on $\Ks$, and a unique $1$-form $\xi\in\Gamma(\TM^{*})$ \st:
\begin{alignat*}{2}
&A=A_{0} 
  &&\quad\mbox{is the unique minimum point of }\B^{A}\phi\\
&A=A_{0}-3i\,\xi
  &&\quad\mbox{is the unique minimum point of }\Sym\B^{A}\phi\\
&A=A_{0}+5i\,\xi 
  &&\quad\mbox{is the unique minimum point of }\Alt\B^{A}\phi\\
&A=A_{0}+15i\,\xi 
  &&\quad\mbox{is the unique vanishing point of }\D^{A}\phi=\tr(\B^{A}\phi) 
\end{alignat*}
All these minimizing connections lie on a same affine line in $\Con(\Ks)$. \\
If any two of them happen to coincide, then all of them must coincide.
\end{theorem}
\begin{theorem}
\label{thm-exact}
The formula for $\xi$ above is
\[ \xi=-\,\tfrac{4}{15\norm{\phi}^{2}}\,d^{*}\sigma(\phi) \]
\end{theorem}
Note that $d^{*}\sigma(\phi)=*\,d\,\sigma(\phi)$, since $\sigma(\phi)$ 
is \SD.

The minima above are determined point-wise.
For a fixed $\phi$, let ${\cal T}$ denote any one of 
$\B\phi$, $\Alt\B\phi$, $\Sym\B\phi$, $\Symo\B\phi$, or 
$g\tens\rec{4}\tr\B\phi$.
Then the form ${\cal T}^{A}\rest{x}$ is an element of 
$\TM^{*}\tens \TM^{*}\tens \TM\rest{x}$.
The latter has a natural \emph{real} inner product induced from $g$, 
given by 
$\Inner{T_{1},T_{2}}
  =\sum_{j,k=1}^{4} \Inner{T_{1}(e_{j},e_{k}),\ T_{2}(e_{j},e_{k}) }$
for any $g$-orthonormal basis 
$\ebasis$ in $\TM\rest{x}$.
(All inner products that appear in this paper are \emph{real}-valued.)

On the other hand, 
the space $\Con(\Ks)$ of all connections $A$ on $\Ks$ is affine, with model 
space $i\,\Gamma(\TM^{*})$. 
The map $A\maps{\cal T}^{A}\rest{x}$ is an affine map
$\Con(\Ks) \to \TM^{*}\tens \TM^{*}\tens\TM\rest{x}$.
The image of $\Con(\Ks)$ through this map is an affine subspace ${\cal C}_{x}$ in $\TM^{*}\tens \TM^{*}\tens\TM\rest{x}$.
Thus, the minimum of the map is \emph{unique}, 
and is exactly the point of ${\cal C}_{x}$ that is closest to the origin of $\TM^{*}\tens \TM^{*}\tens\TM\rest{x}$.
It is the point where the perpendicular from the origin hits ${\cal C}_{x}$, see Figure~\ref{fig-minimum}.
\begin{figure}
\begin{picture}(275,155)(-50,-5)

\thicklines
\put(0,100){\line(2,1){100}}
\put(0,100){\line(2,-1){210}}
\thinlines
\put(50,0){\makebox(0,0){$\bullet$}}
\put(50,0){\line(1,2){50}}
\put(100,100){\makebox(0,0){$\bullet$}}
\put(100,100){\line(2,-1){58}}
\put(160,69.5){\makebox(0,0){$\circ$}}
\put(95,90){\line(2,-1){10}}
\put(110,95){\line(-1,-2){5}}

\put(43,0){\makebox(0,0){0}}
\put(25,100){\makebox(0,0){${\cal C}_x$}}
\put(-50,30){$\TM^{*}\tens \TM^{*}\tens \TM\rest{x}$}
\put(100,109){\makebox(0,0){$min$}}

\end{picture}
\caption{Finding the minimum of 
	$A\maps{\cal T}^{A}\rest{x}$.}
\label{fig-minimum}
\end{figure}
Therefore the map $A\maps{\cal T}^{A}\rest{x}$ attains an 
absolute minimum at $A_{0}$ \Iff, for any other connection $A$, we have:
$\Inner{{\cal T}^{A_{0}}\phi,\  {\cal T}^{A_{0}}\phi- {\cal T}^{A}\phi}\srest{x}
  = 0$. 
Or, by decoding the inner product, if we have:
\begin{equation}
\label{eq-minima}
\tsum_{j,k=1}^{4}
  \Inner{({\cal T}^{A_{0}}\phi)\Sub{e_{j}}e_{k},\ 
  ({\cal T}^{A_{0}}\phi)\Sub{e_{j}}e_{k}
  - ({\cal T}^{A}\phi)\Sub{e_{j}}e_{k}}\srest{x}
  =0
\end{equation}
for some $g$-orthonormal basis 
$\ebasis$ in $\TM\rest{x}$.

Applying the above and using the formulae \eqref{eq-change} 
yields, after some elementary computations:

\begin{lemma}%
\label{lemma-minima}
Let $(B^{A}\phi)\Sub{X}Y=Y\cli\bnabla^{A}_{X}\phi$ and vary $A$. 
Then:
\medskip\\ 
{\bf A.}
The comparison form $\B^{A_{0}}\phi\rest{x}$ is minimal 
in $\TM^{*}\tens\TM^{*}\tens\TM\rest{x}$
\Iff, for all $X\in\TM\rest{x}$, we have:
\[ \Inner{\bnabla^{A_{0}}_{X}\phi,\ i\phi}\srest{x}=0 \]
{\bf B.}
The form $\Sym\B^{A_{0}}\phi\rest{x}$ is minimal 
in $\TM^{*}\tens\TM^{*}\tens\TM\rest{x}$
\Iff, for all $X\in \TM\rest{x}$, we have:
\[ \Inner{\D^{A_{0}}\phi,\ i\phi(X)}\srest{x}
  =6\Inner{\bnabla^{A_{0}}_{X}\phi,\ i\phi}\srest{x} \]  
{\bf C.}
The form $\Alt\B^{A_{0}}\phi\rest{x}$ is minimal 
in $\TM^{*}\tens\TM^{*}\tens\TM\rest{x}$
\Iff, for all $X\in \TM\rest{x}$, we have
\[ \Inner{\D^{A_{0}}\phi,\ i\phi(X)}\srest{x}
  =-2\Inner{\bnabla^{A_{0}}_{X}\phi,\ i\phi}\srest{x} \]
{\bf D.}
Finally, 
$g\tens\rec{4}\tr(\B^{A_{0}}\phi\,\rest{x})
  =g\tens\rec{4}\D^{A_{0}}\phi\,\rest{x}$ 
is minimal in $\TM^{*}\tens\TM^{*}\tens\TM\rest{x}$ 
\Iff\ 
\[ \D^{A_{0}}\phi\,\rest{x}=0 \]
\end{lemma}
Of course, \ref{lemma-minima}.D is utterly trivial. It is included 
only for completeness.

\begin{remark}
\label{rk-min.global}
All the above point-wise minima can be realized globally on $M$ 
if $\phi$ is nowhere-zero. 
For example, for minimizing $\B\phi$, start with a random connection 
$A$ and define 
$\theta_{0}(X)=-\,\rec{\norm{\phi}^{2}}\,\Inner{\bnabla^{A}_{X}\phi,\,i\phi}$.
Then $A_{0}=A+2i\theta_{0}$ is the unique connection that minimizes $\B\phi$.
If $\phi$ has zeros, though, then the globally 
minimizing connection $A$ might explode at the zeros.
(Similarly for $\Sym\B\phi$, $\Alt\B\phi$ and $\D\phi$.)
Nonetheless, for simplicity 
from now on we will talk only of global minimizing connections,
implicitly restricting away from the possible singularities.
\end{remark}
\begin{Proof}{\Corref{cor-more.symplectic}}
\Lemmaref{lemma-minima}.A gives meaning to the 
condition $\inner{\bnabla\phi,\,i\phi}=0$ from \Propref{prop-old}.B: 
it insures the minimality of $\B^{A}\phi$. Together with 
the condition $\D^{A}\phi=0$, 
it implies that all the minimizing connections from 
\Thmref{thm-minima.line} coincide (see also \ref{thm-exact}). 
We could thus rephrase \Propref{prop-old}.B as:
\smallskip

\noindent
\emph{%
The form $\phi^{*}\omega$ is symplectic \Iff\ there is a connection 
$A$ that simultaneously minimizes $\B^{A}\phi$, $\Alt^{A}\B\phi$, 
$\Sym^{A}\B\phi$, and $\D^{A}\phi$.
}\smallskip

In particular, it is equivalent to:
\smallskip

\noindent
\emph{%
The form $\phi^{*}\omega$ is symplectic \Iff\ there is a connection 
$A$ such that $\B^{A}\phi$ and $\Alt\B^{A}\phi$
are both minimal.
}\smallskip

That can be stated:
\smallskip

\noindent
\emph{%
The form $\phi^{*}\omega$ is symplectic 
\Iff\ there is a connection $\tnabla^{A}$ on $\TM=\W^{-}$
which, \emph{through $\phi$}, is simultaneously a closest match to $\nabla$ 
and has the torsion minimized.
}\smallskip

In the special case when $\phi=id:\TM\to\TM$, we read: 
\smallskip

\noindent
\emph{%
The manifold $(M,g,\omega)$ is almost-\kah\ \Iff\ there is a 
$\Ce$-linear $g$-metric connection $\tnabla$ on $\TM$ with 
$\tnabla\rest{\Lambda^{-}}=\nabla\rest{\Lambda^{-}}$,
which is simultaneously closest to $\nabla$ and has minimal torsion.
}\smallskip

This last statement is exactly \Corref{cor-more.symplectic}.
\end{Proof}

\begin{Proof}{\Thmref{thm-minima.line}}
We relate the minimizing connections through the conditions from 
\Lemmaref{lemma-minima}.
Assume first that $A$ is the minimizing connection for $\B\phi$, 
\ie assume that, for all $X$,
\[ \Inner{\bnabla^{A}_{X}\phi,\ i\phi} =0 \]
Assume also that $A+2i\theta$ is the minimizing connection for $\Alt\B\phi$, 
or that, for all $X$,
\[ 2\,\Inner{\bnabla^{A+2i\theta}_{X}\phi,\ i\phi}
  +\Inner{\D^{A+2i\theta}\phi,\ i\phi(X)} =0 \]
Since $\bnabla^{A+2i\theta}_{X}\phi=\bnabla^{A}_{X}\phi+i\theta(X)\phi$, 
we get
\[ 2\norm{\phi}^{2}\theta(X)
  +\Inner{\D^{A+2i\theta}\phi,\ i\phi(X)} =0 \]
Identifying $\TM^{*}$ and $\TM$, we write $\theta(X)=\inner{\theta,X}$, 
and then we have $\norm{\phi}^{2}\inner{\theta,X}=\Inner{i\phi(\theta),\ i\phi(X)}$. 
Therefore, for all $X$:
\begin{gather*}
2\Inner{i\theta\cli\phi,\ i\phi(X)}+
   \Inner{\D^{A+2i\theta}\phi,\ i\phi(X)} =0 \\
\Inner{\D^{A+6i\theta}\phi, \ i\phi(X)} =0  
\end{gather*}
since $\D^{A+6i\theta}\phi=\D^{A+2i\theta}\phi+2i\theta\cli\phi$.
But then, on the support of $\phi$, we must have:
\[ \D^{A+6i\theta}\phi =0 \]
and so we proved:

($\alpha$) \emph{%
If $A$ is the minimum point for $\B\phi$, and $A+2i\theta$ is the 
minimum point for $\Alt\B\phi$, then $\D^{A+6i\theta}\phi=0$.
}

Assume now that $A$ is the minimizing connection for $\B\phi$, 
so that, for all $X$,
\[ \Inner{\bnabla^{A}_{X}\phi,\ i\phi}=0 \]
Assume also that $A+2i\eta$ minimizes $\Sym\B\phi$: for all $X$
\begin{align*}
6\,\Inner{\bnabla^{A+2i\eta}_{X}\phi,\ i\phi}
  -\Inner{\D^{A+2i\eta}\phi,\ i\phi(X)} &=0\\
6\norm{\phi}^{2}\eta(X)
  -\Inner{\D^{A+2i\eta}\phi,\ i\phi(X)} &=0\\
\Inner{6i \eta\cli\phi,\ i\phi(X)}
  -\Inner{\D^{A+2i\eta},\ i\phi(X)} &=0\\
-\Inner{\D^{A-10i\eta}\phi,\ i\phi(X)} &=0
\end{align*}
and so $\D^{A-10i\eta}\phi=0$ on the support of $\phi$.
We obtained:

($\beta$) \emph{%
If $A$ is the minimum point for $\B\phi$, and $A+2i\eta$ is the 
minimum point for $\Sym\B\phi$, then we have $\D^{A-10i\eta}\phi=0$.
}

If we combine now ($\alpha$) and ($\beta$) above with \Rkref{rk-min.global},
then \Thmref{thm-minima.line}.A will follow.
\end{Proof}

\begin{Proof}{\Thmref{thm-exact}}
The formula for $\xi$ is obtained by comparing the minimizing 
connection for $\B\phi$ with the vanishing connection for $\D\phi$.
The main ingredient is formula \eqref{eq-D.formula} (from the 
Introduction; it is proved in \cite{scorpan.harmonic}), combined with
\Lemmaref{lemma-minima}.A and the suitable formula from 
\eqref{eq-change}. Concretely, if $\B^{A}\phi$ is minimal and 
$\D^{A+2i\theta}\phi=0$, then
\[ \theta=-\,\tfrac{2}{\norm{\phi}^{2}}\,d^{*}\sigma \]
Fitting $\theta$ to $\xi$ from \Thmref{thm-minima.line} yields the result.
\end{Proof}

\begin{Proof}{\Lemmaref{lemma-minima}.A}
The minimality formula \eqref{eq-minima}
applied to  ${\cal T}=\B\phi$ becomes the condition 
that, for all $\theta\in \TM^{*}\rest{x}$, we have:
\begin{align*}
\tsum_{j,k=1}^{4} \Inner{(\B^{A_{0}}\phi)\Sub{e_{j}}e_{k},
  \ i\theta(e_{j})\phi(e_{k})}\srest{x} &=0 \\
\tsum_{j,k=1}^{4} \theta(e_{j})
  \Inner{(\B^{A_{0}}\phi)\Sub{e_{j}}e_{k},\ i\phi(e_{k})}\srest{x} &=0
\end{align*}
Since that must happen for all $\theta\in\TM^{*}\rest{x}$, 
we must have, \emph{for all $j$}:
\begin{align*}
\tsum_{k=1}^{4} 
  \Inner{(\B^{A_{0}}\phi)\Sub{e_{j}}e_{k},\ i\phi(e_{k})}\srest{x} &=0 \\
\tsum_{k=1}^{4} 
  \Inner{e_{k}\cli \hnabla^{A_{0}}_{e_{j}}\phi,
  \ e_{k}\cli i\phi}\srest{x} &=0\\
\tsum_{k=1}^{4}
  \Inner{\hnabla^{A_{0}}_{e_{j}}\phi,\, i\phi}\srest{x} &=0
\end{align*}
which means that 
$\Inner{\hnabla^{A}_{e_{j}}\phi,\,i\phi}=0$ and thus concludes.
\end{Proof}

\begin{Proof}{\Lemmaref{lemma-minima}.B}
The minimality condition for $\Sym\B^{A}\phi$ is:
\begin{align*}
\tsum_{j,k=1}^{4} \Inner{(\Sym\B^{A_{0}}\phi)\Sub{e_{j}}e_{k},
  \ (\Sym\B^{A_{0}}\phi)\Sub{e_{j}}e_{k}
    -(\Sym\B^{A}\phi)\Sub{e_{j}}e_{k}}\srest{x} &=0\\
\intertext{%
for all $\theta\in \TM^{*}\rest{x}$. That is:}
\tsum_{j,k=1}^{4} 
  \Inner{(\B^{A_{0}}\phi)\Sub{e_{j}}e_{k}+(\B^{A_{0}}\phi)\Sub{e_{k}}e_{j},
  \ i\theta(e_{j})\phi(e_{k})+i\theta(e_{k})\phi(e_{j})}\srest{x} &=0\\
\tsum_{j=1}^{4} \theta(e_{j}) 
  \bigl(\tsum_{k=1}^{4}\Inner{(\B^{A_{0}}\phi)\Sub{e_{j}}e_{k},\ i\phi(e_{k})}
  +\Inner{(\B^{A_{0}}\phi)\Sub{e_{k}}e_{j},\ i\phi(e_{k})}\srest{x} \bigr) &=0
\end{align*}
Therefore, for all $j$:
\begin{align*}
\tsum_{k=1}^{4} 
  \Inner{(\B^{A_{0}}\phi)\Sub{e_{j}}e_{k},\ i\phi(e_{k})}
  +\Inner{(\B^{A_{0}}\phi)\Sub{e_{k}}e_{j},\ i\phi(e_{k})}\srest{x} &=0\\
\tsum_{k=1}^{4} \Inner{e_{k}\cli\hnabla^{A_{0}}_{e_{j}}\phi,\ e_{k}\cli i\phi}
  +\Inner{e_{j}\cli\hnabla^{A_{0}}_{e_{k}}\phi,\ e_{k}\cli i\phi}\srest{x} &=0
\end{align*}
But $\inner{e_{j}\cli a,\,e_{k}\cli b}=-\inner{e_{k}\cli a,\,e_{j}\cli b}$ 
when $j\neq k$.  So, for every $j$:
\begin{align*}
5\Inner{\hnabla^{A_{0}}_{e_{j}}\phi,\ i\phi}
  -\tsum_{k\neq j}
  \Inner{e_{k}\cli\hnabla^{A_{0}}_{e_{k}}\phi,\ e_{j}\cli i\phi}\srest{x} &=0\\
6\Inner{\hnabla^{A_{0}}_{e_{j}}\phi,\ i\phi}
  -\Inner{\D^{A_{0}}\phi,\ e_{j}\cli i\phi}\srest{x} &=0
\end{align*}
which concludes.
\end{Proof}

\begin{Proof}{\Lemmaref{lemma-minima}.C}
Minimality of $\Alt\B^{A}\phi$ is insured by:
\begin{align*}
\tsum_{j\neq k} 
  \Inner{(\B^{A_{0}}\phi)\Sub{e_{j}}e_{k}-(\B^{A_{0}}\phi)\Sub{e_{k}}e_{j},
  \ i\theta(e_{j})\phi(e_{k})-i\theta(e_{k})\phi(e_{j})}\srest{x} &=0\\
\tsum_{j\neq k} \theta(e_{j})
  \bigl( \Inner{(\B^{A_{0}}\phi)\Sub{e_{j}}e_{k},\ i\phi(e_{k})}
  -\Inner{(\B^{A_{0}}\phi)\Sub{e_{k}}e_{j},\ i\phi(e_{k})}\srest{x} \bigr) &=0
\end{align*}
Then, for all $j$:
\begin{align*}
\tsum_{k\neq j} \bigl(
   \Inner{e_{k}\cli\hnabla^{A_{0}}_{e_{j}}\phi,\ e_{k}\cli i\phi}
  -\Inner{e_{j}\cli\hnabla^{A_{0}}_{e_{k}}\phi,\ e_{k}\cli i\phi}\srest{x} 
  \bigr) & =0\\
3\Inner{\hnabla^{A_{0}}_{e_{j}}\phi,\ i\phi}
  +\tsum_{k\neq j}\Inner{e_{k}\cli\hnabla^{A_{0}}_{e_{k}}\phi,
  \ e_{j}\cli i\phi}\srest{x} &=0\\
2\Inner{\hnabla^{A_{0}}_{e_{j}}\phi,\ i\phi}
  +\Inner{\D^{A_{0}}\phi,\ e_{j}\cli i\phi}\srest{x} &=0
\end{align*}
which concludes.
\end{Proof}

\begin{Proof}{\Lemmaref{lemma-spin.connection} 
  (and construction of $\bnabla^{A}$ and $\tnabla^{A}$)}
In order to prove \ref{lemma-spin.connection}, we need to clearly 
explain how $\nabla$ and $A$ induce the connections $\bnabla^{A}$ 
and $\tnabla^{A}$. In tune with the rest of the paper, we will use 
the language of quaternions.

The group $SO(4)=\Sph{3}\times\Sph{3}\,\big/\pm1$ acts on 
$\aR^{4}\iso\Ha$ 
by $[\xi_{+},\xi_{-}]\cdot v=\xi_{+}v\,\xi_{-}^{-1}$. Its Lie algebra 
is $\lie{so}(4)=\Im\Ha\oplus\Im\Ha$ and its adjoint action on $\aR^{4}$ 
is $(q_{+}\oplus q_{-})\cdot v=q_{+}v-v\,q_{-}$.
Therefore, any $g$-metric connection $D$ on $\TM$ can be written 
locally as $D v=\del v+b_{+}v-v\,b_{-}$ for suitable 
local $1$-forms $b_{\pm}\in\Gamma(\TM^{*}\tens\Im\Ha)$.
In particular, we write the \LC\ connection $\nabla$ on $\TM$ locally as 
$\nabla v=\del v+\goth{a}_{+}v-v\,\goth{a}_{-}$, for some 
$\goth{a}_{\pm}\in\Gamma_{loc}(\TM^{*}\tens\Im\Ha)$.
The connection $D$ induces connections on all tensor bundles of $M$, 
and in particular on $\Lambda^{\pm}(\TM^{*})$; the latter can be 
written locally as $D\rest{\Lambda^{\pm}}f=\del f+b_{\pm}f-f\,b_{\pm}$.

The group $U(1)=\Sph{1}$ acts on $\Ce$ by left multiplication, and its
Lie algebra is $\lie{u}(1)=i\,\aR$. Thus, any unitary connection on the 
complex-line bundle $\Ks$ can be written locally as 
$Az=\del z+2i\alpha z$, for some suitable 
$\alpha\in\Gamma_{loc}(\TM^{*})$.

The group $U(2)=\Sph{1}\times\Sph{3}\,\big/\pm1$ acts on $\Ce^{2}\iso\Ha$
by $[\lambda,\xi]\cdot w=\lambda w\xi^{-1}$. Its Lie algebra is 
$\lie{u}(2)=i\,\aR\oplus\Im\Ha$ and its adjoint action on $\Ce^{2}$ is
$(i\ell\oplus q)\cdot w=i\ell w-wq$. Hence, any unitary connection 
on $\W^{\pm}$ can be written locally as $D w=\del w+i\beta w-w\delta$ for 
some suitable $\beta\in\Gamma_{loc}(\TM^{*})$ and 
$\delta\in\Gamma_{loc}(\TM^{*}\tens\Im\Ha)$.

The group $\Spinc(4)=\Sph{1}\times\Sph{3}\times\Sph{3}$ has Lie 
algebra $\lie{spin}^{\Ce}(4)=i\,\aR\oplus\Im\Ha\oplus\Im\Ha$. 
The canonical maps $\Spinc(4)\to SO(4)$, $\Spinc(4)\to\Sph{1}$, and 
the two maps $\Spinc(4)\to U(2)$ have the Lie algebra versions
\[\begin{CD}
\lie{u}(1) @<<< \lie{spin}^{\Ce}(4) @>>> \lie{so}(4) 
  &\qquad\quad&
  \lie{spin}^{\Ce}(4) @>>> \lie{u}(2)\\
2i\ell && i\ell\oplus q_{+}\oplus q_{-} && q_{+}\oplus q_{-}
  &&
  i\ell\oplus q_{+}\oplus q_{-} && i\ell\oplus q_{\pm}
\end{CD}\]

The consequence is that, if one chooses a unitary connection 
$Az=\del z+2i\alpha z$ on $\Ks$, 
then, combining with the \LC\ connection $\nabla$ and using the 
above maps, $A$ determines \emph{unique} unitary connections on 
$\W^{\pm}$, namely $\bnabla^{A} w=\del w+i\alpha w-w\,\goth{a}_{+}$ on $\Wp$, 
and $\tnabla^{A} w=\del w+i\alpha w-w\,\goth{a}_{-}$ on $\W^{-}$. 

The fact that $\nabla$ and $\tnabla^{A}$ share the coefficient-form 
$\goth{a}_{-}$ can be expressed as 
$\nabla\rest{\Lambda^{-}}=\tnabla^{A}\rest{\Lambda^{-}}$. Indeed
the bundle 
$\Lambda^{-}(\TM^{*})$ has model-fiber $\Im\Ha$ with structure group 
$SO(3)=\Sph{3}\big/\pm1$ acting as $[\xi]\cdot a=\xi a \xi^{-1}$. The 
Lie algebra $\lie{so}(3)$ is $\Im\Ha$ and its adjoint action is 
$h\cdot a=ha-ah$. The cocycle of $\Lambda^{-}$ is induced from 
the cocycle of $\TM$ via the map $SO(4)\to SO(3)$, 
$[\xi_{+},\xi_{-}]\maps [\xi_{-}]$. 
Thus, a connection $Dv=\del v+b_{+}v+v\,b_{-}$ induces on 
$\Lambda^{-}$ the connection $D\rest{\Lambda^{-}}a=\del a+b_{-}a-a\,b_{-}$.

Assume that $D$ is a $g$-metric connection on $\TM$ that 
has $D\rest{\Lambda^{-}}=\nabla\rest{\Lambda^{-}}$. 
Then $Dv=\del v+b_{+}v-v\,\goth{a}_{-}$. If further $D$ is 
$\Ce$-linear, then the quaternionic-imaginary form $b_{+}$ must 
in fact be just complex-imaginary, and thus $b_{-}=i\alpha$ for some 
$\alpha\in\Gamma_{loc}(\TM^{*})$. But it can be verified that
these forms $\alpha$ define a unitary connection $Az=\del z+2i\alpha z$ 
on $\Ks$ (compare with \Rkref{rk-find.A}).
\end{Proof}

\subsection*{Acknowledgments}
We wish to thank R.~Kirby for help and support during sometimes difficult moments. We also wish to thank D.~Freed and C.~Taubes for sharp comments on an earlier version of this paper.



\providecommand{\bysame}{\leavevmode\hbox to3em{\hrulefill}\thinspace}
\providecommand{\MR}{\relax\ifhmode\unskip\space\fi MR }
\providecommand{\MRhref}[2]{%
  \href{http://www.ams.org/mathscinet-getitem?mr=#1}{#2}
}
\providecommand{\href}[2]{#2}

\end{document}